\documentclass[preprint,3p,times]{elsarticle}

\usepackage[T1]{fontenc}
\usepackage{lmodern}
\usepackage{amsmath,amssymb,amsthm,mathtools}
\usepackage{booktabs}
\usepackage{needspace}
\usepackage{array}
\usepackage{enumitem}
\usepackage{microtype}
\usepackage{aliascnt}
\usepackage{xurl}
\usepackage[hidelinks]{hyperref}
\usepackage[capitalize,noabbrev]{cleveref}

\journal{Fuzzy Sets and Systems}
\biboptions{numbers,sort&compress}

\newtheorem{theorem}{Theorem}[section]
\newaliascnt{proposition}{theorem}
\newtheorem{proposition}[proposition]{Proposition}
\aliascntresetthe{proposition}
\newaliascnt{lemma}{theorem}
\newtheorem{lemma}[lemma]{Lemma}
\aliascntresetthe{lemma}
\newaliascnt{corollary}{theorem}
\newtheorem{corollary}[corollary]{Corollary}
\aliascntresetthe{corollary}
\theoremstyle{definition}
\newaliascnt{definition}{theorem}
\newtheorem{definition}[definition]{Definition}
\aliascntresetthe{definition}
\newaliascnt{example}{theorem}
\newtheorem{example}[example]{Example}
\aliascntresetthe{example}
\theoremstyle{remark}
\newaliascnt{remark}{theorem}
\newtheorem{remark}[remark]{Remark}
\aliascntresetthe{remark}

\crefname{theorem}{Theorem}{Theorems}
\Crefname{theorem}{Theorem}{Theorems}
\crefname{proposition}{Proposition}{Propositions}
\Crefname{proposition}{Proposition}{Propositions}
\crefname{lemma}{Lemma}{Lemmas}
\Crefname{lemma}{Lemma}{Lemmas}
\crefname{corollary}{Corollary}{Corollaries}
\Crefname{corollary}{Corollary}{Corollaries}
\crefname{definition}{Definition}{Definitions}
\Crefname{definition}{Definition}{Definitions}
\crefname{example}{Example}{Examples}
\Crefname{example}{Example}{Examples}
\crefname{remark}{Remark}{Remarks}
\Crefname{remark}{Remark}{Remarks}

\newcommand{\I}{\mathcal{I}}
\newcommand{\Psh}{\mathcal{P}}
\newcommand{\FS}{\mathsf{FS}}

\newcommand{\CONT}{\mathsf{CONT}}
\newcommand{\YC}{\mathsf{YC}}
\newcommand{\id}{\mathrm{id}}
\newcommand{\supop}{\operatorname{sup}}

\newcommand{\fs}{\operatorname{fs}}
\newcommand{\Finp}{\operatorname{Fin}^{+}}
\newcommand{\Ylim}{\operatorname{Ylim}}
\newcommand{\ev}{\operatorname{ev}}

\newcommand{\Down}{\mathord{\Downarrow}}

\newcommand{\subdeg}{\operatorname{sub}}
\newcommand{\tn}{\mathbin{\&}}

\setlist[enumerate]{leftmargin=*,itemsep=2pt,topsep=4pt,parsep=0pt}
\allowdisplaybreaks[1]
\numberwithin{equation}{section}

\begin{document}
\begin{frontmatter}
\title{Quantitative FS-domains over continuous \texorpdfstring{\(t\)}{t}-norms\tnoteref{t1}}
\tnotetext[t1]{Research supported by NSF of China (No. 12601126).}

\author[addr1]{Yuxu Chen}
\ead{chenyuxu@scu.edu.cn}
\author[addr1]{Di Yang\texorpdfstring{\corref{cor}}{*}}
\ead{diyang@stu.scu.edu.cn}
\cortext[cor]{Corresponding author}
\address[addr1]{School of Mathematics, Sichuan University, Chengdu 610064, China}

\begin{abstract}
We introduce quantitative FS-domains over continuous \(t\)-norms using a scalar-valued finite-separation degree and an approximate identity. The construction recovers classical FS-domains over the Boolean base. Symmetric quantitative FS-domains are characterized by compactness of their ball topology. In particular,  a complete metric space is a quantitative FS-domain if and only if it is compact. 
Quantitative FS-domains are continuous and form a symmetric monoidal closed category under tensor products and function spaces of Yoneda-continuous maps. This category is Cartesian closed exactly for the Gödel \(t\)-norm.
\end{abstract}
\end{frontmatter}

\section{Introduction}\label{sec:introduction}

A central problem in domain theory is to combine continuous approximation with well-behaved function spaces. Although dcpos and Scott-continuous maps form a Cartesian closed category, the full subcategory of continuous dcpos does not. Classical FS-domains reconcile these requirements through finitely separated approximations of the identity \cite{AbramskyJung1994,Jung1990}.

Finite separation also connects domain theory with topology: classical FS-domains are coherent and Lawson compact \cite[Theorem~4.2.18 and Proposition~4.2.20]{AbramskyJung1994}, and admit descriptions through uniform approximation \cite{JungSunderhauf1998,Kummetz1999}. Quantale-enriched approximation places order-theoretic and metric notions of continuity in a common setting \cite{HofmannWaszkiewicz2011}. This suggests seeking an intrinsic quantitative finite-approximation condition that reflects compactness in the symmetric setting and is stable under enriched function spaces.

Using frame-valued directed ideals, Liu and Zhao established Cartesian closedness of the categories of fuzzy continuous lattices and fuzzy algebraic lattices with fuzzy Scott continuous maps \cite[Theorems~4.9 and~4.12]{LiuZhao2014}. We work with separated fuzzy orders over a continuous \(t\)-norm \(T\), using forward-Cauchy ideals and Yoneda completeness. Existing results establish symmetric monoidal closedness for Yoneda-complete objects \cite{LaiZhang2016} and identify Cartesian closed subcategories through restrictions on hom-values \cite{LaiLuo2025}. Continuous objects in the forward-Cauchy framework, however, are not closed under the standard function spaces. We therefore seek an intrinsic finite-approximation condition that combines continuous approximation with closure under tensor products and Yoneda-continuous function spaces.

We introduce quantitative FS-domains using a scalar-valued finite-separation degree and an approximate identity. For Boolean-valued objects, the definition recovers classical FS-domains. For separated symmetric objects, quantitative FS is equivalent to compactness of the natural ball topology, for every continuous \(t\)-norm. The product \(t\)-norm gives compact metric spaces as a particular case. Thus the same finite-approximation framework encompasses classical FS approximation and symmetric compactness.

We prove that quantitative FS-domains are continuous, and are closed under Yoneda-continuous function spaces. 
Together with closure under the tensor product induced by \(T\), this makes quantitative FS-domains a symmetric monoidal closed category for every continuous \(t\)-norm. 
The function-space theorem recovers classical FS-domain closure in the Boolean case and yields compact function spaces in the symmetric case. Finally, we show that the Cartesian closedness of the category of quantitative FS-domains holds exactly for the G\"odel minimum.

Section~\ref{sec:preliminaries} recalls the preliminaries. Section~\ref{sec:qfs} introduces quantitative FS-domains and proves continuity and symmetric monoidal closedness. Section~\ref{sec:cartesian} classifies Cartesian closedness. Section~\ref{sec:metric} treats symmetric objects, compactness and its function-space consequence.

\section{Preliminaries}\label{sec:preliminaries}
This section recalls the fuzzy-order notions used throughout the paper.  See \cite{LaiZhangZhang2020,Zhang2024} for further background.

\subsection{Fuzzy orders, fuzzy lower sets and suprema}
A continuous \(t\)-norm \(T\) is a continuous, associative, commutative and monotone operation on \([0,1]\) with unit \(1\). We write \(a\tn b=T(a,b)\), and \(a^n\) for the corresponding iterated product, with \(a^0=1\). The symbol \(\tn\) is reserved for the scalar \(t\)-norm product, while \(\otimes\) denotes tensor products of fuzzy ordered sets, maps or ideals. The residual of \(T\) is determined by
\[
 a\Rightarrow b=\bigvee\{c\in[0,1]:a\tn c\leq b\},
 \qquad
 a\tn c\leq b\ \Longleftrightarrow\ c\leq a\Rightarrow b.
\]

A \emph{fuzzy ordered set} with respect to \(T\) \cite[Section~2]{LaiZhangZhang2020} is a nonempty set \(X\) equipped with a fuzzy relation \(X(-,-)\colon X\times X\to[0,1]\) satisfying
\[
 X(x,x)=1,
 \qquad
 X(x,y)\tn X(y,z)\leq X(x,z)
\]
for all \(x,y,z\in X\), together with the separation condition
\[
 X(x,y)=X(y,x)=1\quad\Longrightarrow\quad x=y.
\]
Its underlying order is defined by \(x\leq_0y\) if and only if \(X(x,y)=1\). The fuzzy ordered set \(X\) is \emph{pointed} if this order has a least element, denoted by \(\bot_X\), and it is \emph{symmetric} if \(X(x,y)=X(y,x)\) for all \(x,y\in X\). For the Boolean truth values \(\{0,1\}\), this definition reduces to an ordinary poset.

A map \(f\colon X\to Y\) is \emph{fuzzy-order-preserving} if for all \(x,x'\in X\),
\(
 X(x,x')\leq Y(f(x),f(x'))
\).
Let \([X,Y]\) denote the set of all fuzzy-order-preserving maps from \(X\) to \(Y\). It is a fuzzy ordered set under the pointwise fuzzy order \cite[Section~14]{Zhang2024}
\begin{equation}\label{eq:pointwise-hom}
 [X,Y](f,g)=\bigwedge_{x\in X}Y(f(x),g(x)),
 \qquad \forall f,g\in[X,Y].
\end{equation}
Here \([X,Y](f,g)\) measures the degree to which \(f(x)\) is below \(g(x)\), uniformly in \(x\). In the Boolean case, this is the usual pointwise order on monotone maps. We write \(f\leq_0g\) when \([X,Y](f,g)=1\), and \(\id_X\) for the identity map on \(X\).

A \emph{fuzzy lower set} on \(X\) is a function \(\varphi\colon X\to[0,1]\) satisfying
\(
 X(x,y)\tn\varphi(y)\leq\varphi(x)
\) for all \(x,y\in X\).
Over the Boolean base \(\mathbf2\), fuzzy lower sets are precisely the characteristic functions of ordinary lower sets. 
Let \(\Psh X\) denote the set of all fuzzy lower sets on \(X\). Following the standard notation in  \cite[Section~2]{LaiZhangZhang2020}, define the \emph{fuzzy inclusion degree}
\begin{equation}\label{eq:weight-hom}
 \Psh X(\varphi,\lambda)=\subdeg_X(\varphi,\lambda)
 =\bigwedge_{x\in X}\bigl(\varphi(x)\Rightarrow\lambda(x)\bigr).
\end{equation}
Then \((\Psh X,\subdeg_X)\) is a fuzzy ordered set. For \(r\in[0,1]\),
\begin{equation}\label{eq:bracket-test}
 r\leq\subdeg_X(\varphi,\lambda)
 \quad\Longleftrightarrow\quad \forall x\in X, \
 r\tn\varphi(x)\leq\lambda(x).
\end{equation}
In particular, \(\subdeg_X(\varphi,\lambda)=1\) if and only if \(\varphi\leq\lambda\) pointwise.

Inequalities between functions are understood pointwise. For \(r\in[0,1]\) and \(\varphi\in\Psh X\), define the pointwise scalar multiple \(r\odot\varphi\) by
\(
 (r\odot\varphi)(x)=r\tn\varphi(x).
\)
The scalar multiple \(r\odot\varphi\) is again a fuzzy lower set by \cite[Lemma~2.8]{LaiZhangZhang2020}. Arbitrary joins and meets of fuzzy lower sets are computed pointwise and remain fuzzy lower sets, as follows directly from the defining lower-set inequality. 

For \(x\in X\), the \emph{principal fuzzy lower set} at \(x\) is
\(
 \downarrow_X x=X(-,x).
\)
The principal fuzzy lower sets satisfy 
\begin{equation}\label{eq:yoneda-identities}
 \subdeg_X(\downarrow_X x,\downarrow_X x')=X(x,x'),
 \qquad
 \subdeg_X(\downarrow_X x,\lambda)=\lambda(x),
\end{equation}
and every fuzzy lower set has the pointwise expansion \cite[Corollary~7.11 and Proposition~6.19]{Zhang2024}
\begin{equation}\label{eq:lower-set-expansion}
 \lambda=\bigvee_{x\in X}\lambda(x)\odot\downarrow_X x.
\end{equation}

% Let \(f\colon X\to Y\) be fuzzy-order-preserving. The \emph{direct image} of \(\varphi\in\Psh X\) along \(f\) is
% \begin{equation}\label{eq:direct-image}
%  (f_!\varphi)(y)=\bigvee_{x\in X}\varphi(x)\tn Y(y,f(x)),
% \end{equation}
% while the restriction of \(\lambda\in\Psh Y\) is \(f^*\lambda=\lambda\circ f\). These operations satisfy
% \begin{equation}\label{eq:direct-image-adjunction}
%  \subdeg_Y(f_!\varphi,\lambda)=\subdeg_X(\varphi,f^*\lambda),
% \end{equation}
% as well as
% \begin{equation}\label{eq:direct-image-basic}
%  f_!(\downarrow_X x)=\downarrow_Yf(x),
%  \qquad
%  (gf)_!=g_!f_!.
% \end{equation}

Let \(f\colon X\to Y\) be fuzzy-order-preserving. Following the standard direct-image and inverse-image constructions for fuzzy lower sets
\cite{LaiZhangZhang2020},
we define two maps
\[
f_!\colon\Psh X\to\Psh Y,
\qquad
f^*\colon\Psh Y\to\Psh X.
\]
For \(\varphi\in\Psh X\), its \emph{direct image}
\(f_!\varphi\) is defined by
\begin{equation}\label{eq:direct-image}
(f_!\varphi)(y)
=
\bigvee_{x\in X}
\varphi(x)\tn Y(y,f(x)),
\qquad \forall y\in Y.
\end{equation}
% Thus \(f_!\varphi\) is the fuzzy lower set on \(Y\)
% obtained by transferring \(\varphi\) along \(f\) and
% taking the lower-set closure of its image.
Conversely, for \(\lambda\in\Psh Y\), its \emph{inverse image} \(f^*\lambda\) is defined by
\begin{equation}\label{eq:inverse-image}
f^*\lambda=\lambda\circ f,
\qquad
(f^*\lambda)(x)=\lambda(f(x)),
\quad \forall x\in X.
\end{equation}
Both constructions yield fuzzy lower sets. For \(\varphi\in\Psh X\) and \(\lambda\in\Psh Y\), they satisfy the following identity:
\begin{equation}\label{eq:direct-image-adjunction}
\subdeg_Y(f_!\varphi,\lambda)
=
\subdeg_X(\varphi,f^*\lambda).
\end{equation}
In particular, \(f_!(\downarrow_X x)=\downarrow_Y f(x)\) for every \(x\in X\). Over the Boolean base, \(f_!\) is the lower-set closure of the ordinary direct image, whereas \(f^*\) is the ordinary inverse-image operation.

A point \(c\in X\) is the \emph{supremum} of \(\varphi\in\Psh X\), written \(c=\supop_X\varphi\), if the following equality holds 
\begin{equation}\label{eq:fuzzy-supremum}
X(c,z)=\subdeg_X(\varphi,\downarrow_X z), \qquad\forall z\in X.
\end{equation}
Such a supremum is unique and satisfies \(\varphi\leq\downarrow_X c\). If both \(\supop_X\varphi\) and \(\supop_X\psi\) exist, then
\begin{equation}\label{eq:sup-order-preserving}
 \subdeg_X(\varphi,\psi)
 \leq X(\supop_X\varphi,\supop_X\psi).
\end{equation}
More generally, for a fuzzy-order-preserving map \(F\colon A\to X\) and \(\alpha\in\Psh A\), write \cite[Section~6]{Zhang2024}
\[
 \supop_\alpha F=\supop_X(F_!\alpha)
\]
whenever this supremum exists. Equivalently,
\begin{equation}\label{eq:weighted-supremum}
 X(\supop_\alpha F,z)
 =\subdeg_A\bigl(\alpha,a\mapsto X(F(a),z)\bigr),\qquad \forall z\in X.
\end{equation}

\begin{remark}
In enriched-category terminology, fuzzy ordered sets are separated categories enriched over \(([0,1],\leq,T,1)\), fuzzy lower sets are weights, principal fuzzy lower sets are representable weights, and the supremum in \eqref{eq:weighted-supremum} is a weighted colimit; see \cite{LaiZhang2007,Zhang2024} for the general categorical formulation.
\end{remark}

\subsection{Ideals, completeness and continuous approximation}
A net \((x_i)\) in \(X\) is \emph{forward Cauchy} \cite[Definition~9.1]{Zhang2024} if
\begin{equation}\label{eq:forward-cauchy}
 \bigvee_i\bigwedge_{i\leq j\leq k}X(x_j,x_k)=1.
\end{equation}
A \emph{forward-Cauchy ideal}, or simply an \emph{ideal}, is a fuzzy lower set generated by such a net \cite[Section~3]{LaiZhangZhang2020}:
\begin{equation}\label{eq:generated-weight}
 \varphi(z)=\bigvee_i\bigwedge_{j\geq i}X(z,x_j),\qquad \forall z\in X.
\end{equation}
Let \(\I X\subseteq\Psh X\) denote the set of all forward-Cauchy ideals on \(X\), equipped with the restriction of \(\subdeg_X\). Thus \(\I X\) is itself a fuzzy ordered set.
A \emph{Yoneda limit} of a forward-Cauchy net \((x_i)\), denoted by \(\Ylim_i x_i\), is an element \(c\in X\) satisfying 
\begin{equation}\label{eq:yoneda-limit}
 X(c,z)=\bigvee_i\bigwedge_{j\geq i}X(x_j,z),
 \qquad \forall z\in X.
\end{equation}

\begin{proposition}\cite[Definition~9.9, Theorem~9.10, Lemma~9.11 and Proposition~9.13]{Zhang2024}\label{prop:forward-cauchy-representation}
A fuzzy lower set \(\alpha\) is an ideal if and only if, for all \(r\in[0,1]\) and \(\lambda,\mu\in\Psh X\),
\begin{equation}\label{eq:canonical-laws}
\begin{aligned}
 \subdeg_X(\alpha,r\odot\lambda)
 &=r\tn\subdeg_X(\alpha,\lambda),\\
 \subdeg_X(\alpha,\lambda\vee\mu)
 &=\subdeg_X(\alpha,\lambda)\vee\subdeg_X(\alpha,\mu),\\
 \subdeg_X(\alpha,0)&=0.
\end{aligned}
\end{equation}
Here \(0\) denotes the constant-zero fuzzy lower set. Principal fuzzy lower sets are ideals, direct images preserve ideals, and every ideal is inhabited: \(\bigvee_x\alpha(x)=1\). The supremum of a generated ideal is precisely the Yoneda limit of its generating net, whenever either exists.
\end{proposition}

\begin{lemma}\cite[Proposition~9.13]{Zhang2024}\label{lem:forward-cauchy-bracket}
If \(\varphi\) is generated by a forward-Cauchy net \((x_i)\), then, for every fuzzy lower set \(\lambda\), the real net \((\lambda(x_i))\) converges and
\(\subdeg_X(\varphi,\lambda)=\lim_i\lambda(x_i).
\)
\end{lemma}

A fuzzy ordered set is \emph{Yoneda complete} if every forward-Cauchy net has a Yoneda limit, equivalently, if every ideal has a supremum \cite[Definition~9.3 and Theorem~9.10]{Zhang2024}. A fuzzy-order-preserving map \(f\colon X\to Y\) between Yoneda-complete fuzzy ordered sets is \emph{Yoneda continuous}, or simply \emph{continuous}, if it preserves Yoneda limits of forward-Cauchy nets \cite[Definition~9.3]{Zhang2024}. 
Equivalently \cite[discussion following Proposition~9.13]{Zhang2024},
\begin{equation}\label{eq:ideal-continuity}
 f(\supop_X\varphi)=\supop_Y(f_!\varphi), \qquad \forall
  \varphi\in\I X.
\end{equation}
Unless stated otherwise, completeness of a fuzzy ordered set means Yoneda completeness. 
For a Yoneda-complete \(X\), the \emph{way-below degree} relative to forward-Cauchy ideals \cite[Section~17, before Proposition~17.3]{Zhang2024} is defined by
\begin{equation}\label{eq:relative-way-below}
 w_X(z,x)=\bigwedge_{\varphi\in\I X}
 \bigl(X(x,\supop_X\varphi)\Rightarrow\varphi(z)\bigr).
\end{equation}
The analogous ideal-based formula for fuzzy directed ideals appears in \cite[Definition~3.1]{SuLi2012}; here \(\I X\) always denotes forward-Cauchy ideals.
We call a Yoneda-complete \(X\) a \emph{continuous quantitative domain} \cite[Definition~17.1 and Proposition~17.3]{Zhang2024} if
\begin{equation}\label{eq:continuity-approximation}
 \Down_X x:=w_X(-,x)\in\I X,
 \qquad
 \supop_X(\Down_X x)=x,
 \qquad \forall x\in X.
\end{equation}
Let \(\YC_T\) denote the category of Yoneda-complete fuzzy ordered sets and continuous maps, and let \(\CONT_T\) be its full subcategory of continuous quantitative domains. A subscript \(\bot\) denotes the pointed full subcategory.

\subsection{Products and continuous function spaces}\label{subsec:existing-categories}

The Cartesian product and tensor product of fuzzy ordered sets
\(X\) and \(Y\) have the same carrier \(X\times Y\), with
fuzzy orders defined respectively by
\begin{align}
 (X\times Y)((x,y),(x',y'))
 &=X(x,x')\wedge Y(y,y'),
 \label{eq:cartesian-hom}\\
 (X\otimes Y)((x,y),(x',y'))
 &=X(x,x')\tn Y(y,y').
 \label{eq:tensor-hom}
\end{align}
The former is the categorical product of fuzzy ordered sets.
Indeed, the coordinate projections are fuzzy-order-preserving,
and for any fuzzy-order-preserving maps
\(f\colon Z\to X\) and \(g\colon Z\to Y\), their unique pairing
\[
 \langle f,g\rangle\colon Z\to X\times Y,
 \qquad z\mapsto(f(z),g(z)),
\]
is fuzzy-order-preserving. For Yoneda-complete objects, the same
universal property holds in \(\YC_T\), since ideal suprema in
Cartesian products are computed coordinatewise
(see \Cref{prop:pointwise-completeness}).

In contrast, \(X\otimes Y\) is the monoidal product induced by
\(T\), and need not be a categorical product. The two products
coincide when \(T=\min\), or when either factor is Boolean-valued
\cite{LaiLuo2026,Zhang2024}. For fuzzy-order-preserving maps
\(f\colon X\to X'\) and \(g\colon Y\to Y'\), both
\(f\times g\) and \(f\otimes g\) act coordinatewise on their
respective products.

For Yoneda-complete \(X,Y\), let \([X,Y]_{\I}\) denote the
fuzzy ordered set of all Yoneda-continuous maps from \(X\) to
\(Y\), equipped with the pointwise fuzzy order
\eqref{eq:pointwise-hom}. This set is nonempty. Indeed, for each
\(y\in Y\), the constant map \(c_y\colon X\to Y\) is
fuzzy-order-preserving, and the inhabitedness of ideals gives
\[
 (c_y)_!\varphi=\downarrow_Y y,
 \qquad \forall\varphi\in\I X.
\]
Hence \(c_y\) is Yoneda continuous. For \(x\in X\), write
\[
 \ev_x\colon[X,Y]_{\I}\to Y,
 \qquad \ev_x(h)=h(x).
\]

We next recall the ideal and weighted-supremum formulas needed
below. The external tensor of two forward-Cauchy ideals is
again an ideal, as follows from their net representations
\cite[Proposition~9.13 and the discussion preceding
Lemma~14.6]{Zhang2024}. The interchange formula
\eqref{eq:iterated-fubini} is an instance of the Fubini theorem
for weighted colimits
\cite[Section~3.3, Eq.~(3.21)]{Kelly1982}.
The coordinatewise formula \eqref{eq:sup-coordinatewise} follows
from the representations of ideals by forward-Cauchy nets
\cite[Lemma~14.6 and its proof]{Zhang2024}.

\begin{lemma}\label{lem:external-ideal-suprema}
Let \(\alpha\in\I A\) and \(\beta\in\I B\). Then
\[
 (\alpha\otimes\beta)(a,b)=\alpha(a)\tn\beta(b)
\]
defines an ideal on \(A\otimes B\). If \(Z\) is Yoneda complete and \(F\colon A\otimes B\to Z\) is fuzzy-order-preserving, then
\begin{equation}\label{eq:iterated-fubini}
 \supop_{\alpha\otimes\beta}F
 =\supop_\alpha\bigl(a\mapsto\supop_\beta F(a,-)\bigr)
 =\supop_\beta\bigl(b\mapsto\supop_\alpha F(-,b)\bigr).
\end{equation}
In particular, if \(A\) and \(B\) are Yoneda complete, then
\begin{equation}\label{eq:sup-coordinatewise}
 \supop_{A\otimes B}(\alpha\otimes\beta)
 =\bigl(\supop_A\alpha,\supop_B\beta\bigr).
\end{equation}
\end{lemma}

\begin{proposition}\cite[Section~5]{LaiLuo2026}, \cite[Theorem~14.1 and Lemma~14.2]{Zhang2024}\label{prop:pointwise-completeness}
Finite Cartesian products of Yoneda-complete fuzzy ordered sets are Yoneda complete, with coordinatewise ideal suprema. If \(X,Y\) are Yoneda complete, then \([X,Y]_{\I}\) is Yoneda complete and its ideal suprema are pointwise. The evaluation map, precomposition and postcomposition by continuous maps are continuous.
\end{proposition}

Explicitly, if \(\psi\in\I(X\times Y)\), its coordinate ideals are
\[
 \psi_X=(\pi_X)_!\psi,\qquad \psi_Y=(\pi_Y)_!\psi,
\]
where \(\pi_X(x,y)=x\) and \(\pi_Y(x,y)=y\). They satisfy
\[
 \psi_X(x)=\bigvee_{y\in Y}\psi(x,y),\qquad
 \psi_Y(y)=\bigvee_{x\in X}\psi(x,y),
\]
and
\[
 \supop_{X\times Y}\psi=(\supop_X\psi_X,\supop_Y\psi_Y).
\]
For \(H=[X,Y]_{\I}\) and \(\rho\in\I H\), pointwise computation means
\[
 (\supop_H\rho)(x)=\supop_Y((\ev_x)_!\rho),
 \qquad \forall x\in X.
\]

\begin{theorem}\cite{LaiZhang2016}, \cite[Theorems~14.1 and~14.5]{Zhang2024}\label{thm:ambient-real}
The category \(\YC_T\) is symmetric monoidal closed with tensor \eqref{eq:tensor-hom}, singleton unit and internal hom \([X,Y]_{\I}\). Thus there is a natural bijection
\[
 \YC_T(Z\otimes X,Y)\cong\YC_T(Z,[X,Y]_{\I}),
 \qquad
 \widehat b(z)(x)=b(z,x).
\]
Moreover, \(\YC_T\) is Cartesian closed if and only if \(T=\min\).
\end{theorem}
In particular, tensor products of continuous maps are continuous. 
The completeness and interchange formulas in \Cref{lem:external-ideal-suprema,prop:pointwise-completeness} will be used in the finite-separation arguments.

The distinction between the two closed structures is visible in evaluation. For the pointwise internal hom,
\[
 [X,Y]_{\I}(f,g)\tn X(x,x')\leq Y(f(x),g(x')).
\]
A Cartesian exponential \(Y^X\), if it exists, instead requires
\[
 Y^X(f,g)\wedge X(x,x')\leq Y(f(x),g(x')).
\]
If \(X\) is Boolean-valued, then \(Z\otimes X=Z\times X\) for every \(Z\), so \([X,Y]_{\I}\) is also a Cartesian exponential in \(\YC_T\).

\subsection{The obstruction for continuous quantitative domains}
We now verify that the constructions above recover the classical Boolean case and then transfer a classical function-space obstruction to \(\CONT_T\).

\begin{proposition}\cite[Sections~9 and~17]{Zhang2024}\label{prop:boolean-correspondence}
For a nonempty poset \(P\), put \(P_{\mathbf2}(x,y)=1\) when \(x\leq y\), and \(P_{\mathbf2}(x,y)=0\) otherwise. Its ideals are the characteristic functions of nonempty directed lower sets. Thus \(P_{\mathbf2}\) is Yoneda complete (respectively continuous) exactly when \(P\) is a dcpo (respectively a continuous dcpo), and continuous maps between Boolean-valued objects are precisely Scott-continuous maps. For dcpos \(D,E\),
\begin{equation}\label{eq:boolean-function}
 [D_{\mathbf2},E_{\mathbf2}]_{\I}=[D\to E]_{\mathbf2},
\end{equation}
where \([D\to E]\) is ordered pointwise.
\end{proposition}

\begin{example}\label{ex:complete-not-continuous}
Completeness alone does not imply continuous approximation. Let \(D\) be the dcpo consisting of two increasing \(\omega\)-chains with a common top and no comparisons between the two chains \cite[Figure~4]{AbramskyJung1994}. The top element is not the supremum of its way-below elements. Hence \(D_{\mathbf2}\in\YC_T\setminus\CONT_T\).
\end{example}

Consider the pointed algebraic dcpos
\[
 D_\infty=\{\bot,a,b\}\cup\{c_n:n\in\mathbb N\},
 \qquad
 E_6=\{0,\alpha,\beta,\gamma,\delta,1\},
\]
with orders generated by
\[
 \bot<a,\quad\bot<b,\quad a<c_n,\quad b<c_n,\quad \forall n\in\mathbb N,
 \qquad
 0<\alpha,\beta<\gamma,\delta<1.
\]
The posets \(E_6\) and \(D_\infty\) are the classical examples in \cite[Figures~11 and~12, respectively]{AbramskyJung1994}. Their Scott function space \([D_\infty\to E_6]\) is not continuous \cite[Lemma~4.3.1]{AbramskyJung1994}.

\begin{proposition}\label{prop:CONT-not-CCC}
For every continuous \(t\)-norm \(T\), neither \(\CONT_T\) nor \(\CONT_{T,\bot}\) is Cartesian closed. Neither category is symmetric monoidal closed with respect to \(\otimes\).
\end{proposition}

\begin{proof}
By \Cref{prop:boolean-correspondence},
\((D_\infty)_{\mathbf2},(E_6)_{\mathbf2}
\in\CONT_{T,\bot}\), whereas
\([D_\infty\to E_6]_{\mathbf2}\notin\CONT_T\).

Suppose that a Cartesian exponential
\(H=((E_6)_{\mathbf2})^{(D_\infty)_{\mathbf2}}\)
exists in either category. By the exponential
universal property, its underlying set is
\([D_\infty\to E_6]\).

For every \(x\in D_\infty\), the evaluation
map \(\ev_x:H\to(E_6)_{\mathbf2}\),
\(h\mapsto h(x)\), is a morphism. Hence
\(H(h,k)\leq(E_6)_{\mathbf2}(h(x),k(x))\).
Therefore, \(h\nleq k\) pointwise implies
\(H(h,k)=0\).

Conversely, let \(h\leq k\) pointwise.
Put \(S=\{0<1\}\), and define
\(F:S\times D_\infty\to E_6\) by
\(F(0,x)=h(x)\) and \(F(1,x)=k(x)\).
For every \(U\in\sigma(E_6)\),
\[
F^{-1}(U)
=(S\times h^{-1}(U))
\cup(\{1\}\times k^{-1}(U)).
\]
Thus \(F\) is Scott continuous.
Since \((S\times D_\infty)_{\mathbf2}
=S_{\mathbf2}\times(D_\infty)_{\mathbf2}\)
is the categorical product, \(F\) induces
a morphism into \((E_6)_{\mathbf2}\).
Its transpose \(\widehat F:S_{\mathbf2}\to H\)
satisfies \(\widehat F(0)=h\) and
\(\widehat F(1)=k\). Hence
\(1=S_{\mathbf2}(0,1)\leq H(h,k)\),
so \(H(h,k)=1\).
Consequently,
\[
H(h,k)=
\begin{cases}
1,&h\leq k\text{ pointwise},\\
0,&\text{otherwise}.
\end{cases}
\]
Thus \(H\cong[D_\infty\to E_6]_{\mathbf2}\)
as fuzzy ordered sets, contradicting
\([D_\infty\to E_6]_{\mathbf2}\notin\CONT_T\).
Therefore neither category is Cartesian closed.

Finally, since \((D_\infty)_{\mathbf2}\)
is Boolean-valued, for every \(Z\),
\(Z\otimes(D_\infty)_{\mathbf2}
=Z\times(D_\infty)_{\mathbf2}\).
Hence a monoidal internal hom from
\((D_\infty)_{\mathbf2}\) to
\((E_6)_{\mathbf2}\) would also be a
Cartesian exponential, a contradiction.
Therefore neither category is symmetric
monoidal closed with respect to \(\otimes\).
\end{proof}

The obstruction concerns the full category of all continuous quantitative domains. The quantitative FS-domains introduced below form a symmetric monoidal closed class of continuous objects. A different approach obtains Cartesian closed categories of Yoneda-complete objects by restricting hom-values \cite{LaiLuo2025}.
\section{Quantitative FS-domains}\label{sec:qfs}

This section introduces quantitative FS-domains, derives continuous approximation from their finite-separation approximate identities, and proves closure under tensor products and function spaces of Yoneda-continuous maps, yielding symmetric monoidal closedness. 

Recall that a Scott-continuous endomap \(f\colon D\to D\) on a dcpo \(D\) is \emph{finitely separated from the identity} if there exists a finite nonempty set \(M\subseteq D\) such that, for every \(x\in D\), some \(m\in M\) satisfies
\(
f(x)\leq m\leq x.
\)
Classically, a dcpo is an FS-domain if its identity is the directed pointwise supremum of Scott-continuous endomaps finitely separated from the identity \cite{AbramskyJung1994}. 

\subsection{From finite separating sets to a quantitative degree}

We now replace this Boolean finite-separation condition by a quantitative one. Let \(X\) be a fuzzy ordered set and let \(f\colon X\to X\) be fuzzy-order-preserving. For \(x,m\in X\), the values \(X(f(x),m)\) and \(X(m,x)\) measure the degrees to which
\(
f(x)\leq m
\  \text{and}\
m\leq x
\)
hold. Their combined degree is
\[
X(f(x),m)\tn X(m,x).
\]
For a finite nonempty set \(M\subseteq X\), define
\begin{align}
K_M(z,x)
&=
\bigvee_{m\in M}
X(z,m)\tn X(m,x),
\label{eq:finite-kernel}\\
s_M(f)
&=
\bigwedge_{x\in X}K_M(f(x),x),\\
\fs_X(f)
&=
\bigvee_{M\in\Finp(X)}s_M(f),
\label{eq:finite-separation-degree}
\end{align}
where \(\Finp(X)\) denotes the set of all finite nonempty subsets of \(X\).

For a fixed \(M\), the value \(K_M(f(x),x)\) selects the best intermediate point \(m\in M\) for \(x\), while \(s_M(f)\) takes the worst value over all \(x\in X\). Finally, \(\fs_X(f)\) allows the finite set \(M\) to vary. Thus \(\fs_X(f)\) measures the degree to which \(f\) can be uniformly separated from the identity by finitely many intermediate points.
Over the Boolean base, \(\fs_X(f)=1\) exactly when \(f\) is finitely separated from the identity in the classical sense.

\begin{lemma}\label{lem:fs-basic}
For fuzzy-order-preserving endomaps \(f,g\colon X\to X\),
\begin{equation}\label{eq:fs-lower}
  [X,X](f,g)\tn\fs_X(g)\leq\fs_X(f)
  \leq[X,X](f,\id_X).
\end{equation}
If \(f\leq_0\id_X\) and \(f(X)\) is finite, then \(\fs_X(f)=1\).
\end{lemma}

\begin{proof}
For \(M\in\Finp(X)\), transitivity gives
\[
 K_M(z,x)=\bigvee_{m\in M}X(z,m)\tn X(m,x)\leq X(z,x).
\]
Hence \(s_M(f)\leq[X,X](f,\id_X)\), and 
\(\fs_X(f)=\bigvee_{M\in\Finp(X)}s_M(f)
\leq[X,X](f,\id_X)\).
% taking the join over \(M\) yields the upper bound for \(\fs_X(f)\).

Put \(a=[X,X](f,g)\). For \(x\in X\) and \(m\in M\),
\begin{align*}
 a\tn X(g(x),m)\tn X(m,x) &\leq X(f(x),g(x))\tn X(g(x),m)\tn X(m,x)\\
                          &\leq X(f(x),m)\tn X(m,x).
\end{align*}
Preservation of joins by \(a\tn(-)\) therefore gives
\[
 a\tn s_M(g)\leq a\tn K_M(g(x),x)\leq K_M(f(x),x),
 \qquad \forall x\in X.
\]
Taking the meet over \(x\) yields \(a\tn s_M(g)\leq s_M(f)\). Consequently,
\[
 [X,X](f,g)\tn\fs_X(g)
 =\bigvee_{M\in\Finp(X)}a\tn s_M(g)
 \leq\bigvee_{M\in\Finp(X)}s_M(f)=\fs_X(f).
\]

Finally, if \(f\leq_0\id_X\) has finite image, set \(M=f(X)\). For every \(x\in X\),
\[
 K_M(f(x),x)\geq X(f(x),f(x))\tn X(f(x),x)=1.
\]
Thus \(s_M(f)=1\) and \(\fs_X(f)=1\).
\end{proof}

\subsection{Finite-separation approximate identities}
We next formulate the analogue of the approximate identity used for classical FS-domains. For a Yoneda-complete fuzzy ordered set \(X\), the space
of continuous endomaps
\(
E_X=[X,X]_{\I}
\)
is Yoneda complete by \Cref{prop:pointwise-completeness}.
Here \(\I E_X\) denotes the fuzzy ordered set of ideals on \(E_X\).
An ideal \(\eta_X\) on \(E_X\) represents a directed quantitative family of endomaps. The equation \(\supop_{E_X}\eta_X=\id_X\) says that this family approximates the identity, while \(\eta_X(f)\leq\fs_X(f)\) requires every coefficient assigned to \(f\) to be justified by its finite-separation degree. This is the direct quantitative counterpart of a directed finite-separation approximate identity.

\begin{definition}\label{def:qfs}
A \emph{finite-separation approximate identity} on a Yoneda-complete fuzzy ordered set \(X\) is an ideal \(\eta_X\in\I E_X\) satisfying
\begin{equation}\label{eq:qfs-approximate-identity}
  \supop_{E_X}\eta_X=\id_X,
  \qquad
  \eta_X\leq\fs_X.
\end{equation}
A \emph{quantitative FS-domain} is a Yoneda-complete fuzzy ordered set admitting such an approximate identity. We write \(\FS_T\) for the full subcategory of \(\YC_T\) on quantitative FS-domains, and \(\FS_{T,\bot}\) for its pointed full subcategory. 
\end{definition}

Since suprema in \(E_X\) are pointwise, the first condition in \eqref{eq:qfs-approximate-identity} is equivalent to
\begin{equation}\label{eq:identity-colimit-pointwise}
  X(x,z)=\subdeg_{E_X}\bigl(\eta_X,f\mapsto X(fx,z)\bigr)
  ,\qquad \forall x,z\in X.
\end{equation}
For continuous \(t\)-norms, \Cref{def:qfs} has a direct net formulation.

\begin{proposition}\label{prop:net-qfs}
Let \(X\) be a Yoneda-complete fuzzy ordered set. Then \(X\) is a quantitative FS-domain if and only if there is a forward-Cauchy net \((f_i)\) in \(E_X\) such that
\begin{equation}\label{eq:net-qfs}
  \Ylim_i f_i=\id_X,
  \qquad
  \lim_i\fs_X(f_i)=1.
\end{equation}
\end{proposition}

\begin{proof}
By \Cref{prop:forward-cauchy-representation}, every ideal
\(\eta\in\I E_X\) is generated by a forward-Cauchy net \((f_i)\), and
\(\supop_{E_X}\eta=\Ylim_i f_i\).
Since the restriction \(\fs_X|_{E_X}\) is a fuzzy lower set by \Cref{lem:fs-basic},
\Cref{lem:forward-cauchy-bracket} gives
\[
 \eta\leq\fs_X
 \ \Longleftrightarrow\
 \subdeg_{E_X}(\eta,\fs_X)=1
 \ \Longleftrightarrow\
 \lim_i\fs_X(f_i)=1.
\]
Thus an ideal satisfies \eqref{eq:qfs-approximate-identity} exactly when
its generating net satisfies \eqref{eq:net-qfs}.
\end{proof}

% \paragraph{Relation with familiar classes}
The definition requires a finite-separation approximate identity in addition to Yoneda completeness; continuity will follow from \Cref{thm:qfs-continuous}. Every nonempty finite fuzzy ordered set is quantitative FS by \Cref{prop:finite-objects}. 
% In particular, every finite fuzzy continuous lattice is quantitative FS. 
% We do not claim that every infinite fuzzy continuous lattice is quantitative FS, since the finite-separation approximate identity is an additional uniform approximation requirement. Likewise, completeness and separability are insufficient in the metric setting: \Cref{cor:compact-metric} shows that, for the product \(t\)-norm, the complete separated symmetric metric spaces that are quantitative FS are exactly the compact ones.

\begin{remark}
\label{prop:boolean-recovery}
Over the Boolean base \(\mathbf2\), \Cref{def:qfs} is equivalent to the finite-separation approximate-identity condition recalled above.

For Boolean-valued \(X\), the fuzzy order on \(E_X\) is Boolean and
\(\fs_X(f)\in\{0,1\}\). By \Cref{prop:boolean-correspondence}, every
\(\eta\in\I E_X\) is the characteristic function of a nonempty directed
lower set \(D\subseteq E_X\). Hence
\[
 \supop_{E_X}\eta=\id_X,\quad\eta\leq\fs_X
 \ \Longleftrightarrow\
 \bigvee D=\id_X,\quad \fs_X(f)=1,\quad \forall f\in D.
\]
Since all homs are Boolean,
\[
 \fs_X(f)=1
 \quad\Longleftrightarrow\quad
 \exists M\in\Finp(X)\ \forall x\in X\ \exists m\in M:
 f(x)\leq m\leq x.
\]
Thus the latter condition is precisely a classical finite-separation
approximate identity. Conversely, the lower closure of any such directed
family is directed, has the same supremum, and consists of finitely
separated maps by \eqref{eq:fs-lower}.
Its characteristic function therefore satisfies
\eqref{eq:qfs-approximate-identity}. 
\end{remark}

\begin{proposition}\label{prop:finite-objects}
Every nonempty finite separated fuzzy ordered set over a continuous \(t\)-norm is complete and quantitative FS. Every fuzzy-order-preserving map from it to a complete fuzzy ordered set is continuous.
\end{proposition}

\begin{proof}
Let \(X\) be finite and \(\varphi\in\I X\). By
\eqref{eq:lower-set-expansion} and the ideal identities,
\[
 1=\subdeg_X(\varphi,\varphi)
  =\bigvee_{a\in X}
    \varphi(a)\tn\subdeg_X(\varphi,\downarrow_X a).
\]
The finite maximum is attained. Since \(a\tn b\leq a\wedge b\), some
\(a\in X\) satisfies
\(\varphi(a)=\subdeg_X(\varphi,\downarrow_X a)=1\). Therefore
\[
 X(z,a)=X(z,a)\tn\varphi(a)\leq\varphi(z)\leq X(z,a),
 \qquad \forall z\in X,
\]
so \(\varphi=\downarrow_X a\) and \(\supop_X\varphi=a\).
For every fuzzy-order-preserving \(f\colon X\to Y\) with \(Y\) Yoneda complete,
\[
 \supop_Y(f_!\varphi)
 =\supop_Y(\downarrow_Y f(a))
 =f(a)=f(\supop_X\varphi).
\]
Thus \(X\) is complete and every such \(f\) is continuous.
Taking \(M=X\) gives \(\fs_X(\id_X)=1\). For
\(\eta=E_X(-,\id_X)\), \Cref{lem:fs-basic} yields
\[
 \eta\in\I E_X,\qquad \supop_{E_X}\eta=\id_X,
 \qquad
 \eta(f)=E_X(f,\id_X)\tn\fs_X(\id_X)\leq\fs_X(f).
\]
Hence \(\eta\) is a finite-separation approximate identity.
\end{proof}

\begin{example}\label{ex:product-dual-fs}
Let \(T\) be the product \(t\)-norm, whose implication is
\[
 u\Rightarrow v=
 \begin{cases}
  1,&u\leq v,\\
  v/u,&u>v.
 \end{cases}
\]
Equip \(X=[0,1]\) with \(X(u,v)=v\Rightarrow u\), the dual of the canonical truth-value order \cite[Example~4.3(i)]{Zhang2024}. Its underlying order is the reverse of the usual order, with least element \(1\); it is infinite, non-symmetric and non-Boolean. By \cite[Lemma~9.5]{Zhang2024}, \(X\) is Yoneda complete, and every forward-Cauchy net converges in the usual real topology to its Yoneda limit.

We first compute the finite-separation degree of every \(h\in E_X\):
\[
 \fs_X(h)=
 \begin{cases}
  h(1),&h(0)>0,\\
  0,&h(0)=0.
 \end{cases}
\]
Put \(a=h(0)\) and \(b=h(1)\). Fuzzy order preservation gives \(a\leq h(u)\) and \(bu\leq h(u)\), so \(E_X(h,\id_X)=b\). Hence \(\fs_X(h)\leq b\) by \Cref{lem:fs-basic}.
If \(a=0\), Yoneda continuity gives \(h(u)\to0\) as \(u\downarrow0\). For a finite nonempty subset \(M\subseteq X\) with a positive element, put \(c=\min(M\setminus\{0\})\). For \(u>0\), the term corresponding to \(m=0\), if present, is zero, and every \(m\in M\setminus\{0\}\) satisfies
\[
 X(h(u),m)\tn X(m,u)\leq h(u)/m\leq h(u)/c.
\]
Thus \(s_M(h)=0\); the same holds for \(M=\{0\}\), proving \(\fs_X(h)=0\).

If \(a>0\), let \(M_k=\{a+j(1-a)/k:j=0,\ldots,k\}\) and \(\delta_k=(1-a)/k\). For each \(u\), choose \(m\in M_k\) with \(h(u)\leq m\leq h(u)+\delta_k\). Since \(h(u)\geq a\) and \(u\Rightarrow h(u)\geq b\),
\[
 X(h(u),m)\tn X(m,u)
 =\frac{h(u)}{m}\,(u\Rightarrow m)
 \geq\frac{ab}{a+\delta_k}.
\]
Consequently, \(s_{M_k}(h)\geq ab/(a+\delta_k)\). Letting \(k\to\infty\) proves \(\fs_X(h)=b\), completing the formula.

Now define
\[
 g_n(u)=\frac{nu+1}{n+1},\qquad n\geq1.
\]
For \(0\leq u<v\), the inequality \(u/v\leq(nu+1)/(nv+1)\) shows that \(g_n\) is fuzzy-order-preserving. Its usual continuity and the preceding description of Yoneda limits give \(g_n\in E_X\). Moreover, \(E_X(g_n,g_m)=1\) for \(n\leq m\), so \((g_n)\) is forward Cauchy. Since \(g_n(u)\downarrow u\), its pointwise Yoneda limit is \(\id_X\); this is also its Yoneda limit in \(E_X\), by the pointwise ideal suprema in \Cref{prop:pointwise-completeness}. The formula above gives \(\fs_X(g_n)=1\), because \(g_n(0)>0\) and \(g_n(1)=1\). Thus \(X\) is quantitative FS by \Cref{prop:net-qfs}.

The generated approximate identity is
\[
 \eta(h)=\bigvee_n E_X(h,g_n),\qquad h\in E_X.
\]
It is nonprincipal: \(\eta\leq\fs_X\) by \Cref{lem:fs-basic}, and the formula gives \(\fs_X(\id_X)=0\), hence \(\eta(\id_X)=0\). The principal ideal with supremum \(\id_X\) has value \(1\) there.
\end{example}

\begin{example}\label{ex:product-order-not-fs}
For the product \(t\)-norm, let \(V=[0,1]\) with \(V(u,v)=u\Rightarrow v\). This canonical truth-value order is algebraic, hence continuous \cite[Example~17.8 and Proposition~17.9]{Zhang2024}. For every \(h\in E_V\),
\[
 \fs_V(h)=
 \begin{cases}
  1,&h=c_0,\\
  0,&h\ne c_0,
 \end{cases}
\]
where \(c_0\) is the constant-zero map.

To prove the formula, let \(h\in E_V\) be different from \(c_0\). Then \(h(u)>0\) for every \(u>0\): choose \(a\) with \(h(a)>0\) and use \(V(a,u)>0\) together with fuzzy order preservation. For a finite nonempty subset \(M\subseteq V\) containing positive elements, put \(b=\min(M\setminus\{0\})\). If \(0<u<b\), then
\[
 V(h(u),m)\tn V(m,u)\leq u/m\leq u/b
 \qquad(m\in M\setminus\{0\}),
\]
and the term corresponding to \(m=0\) is \(0\). Thus \(s_M(h)=0\), also when \(M=\{0\}\), and \(\fs_V(h)=0\). On the other hand, \(c_0\) is continuous and \(\fs_V(c_0)=1\) by \Cref{lem:fs-basic}.

If \(V\) were quantitative FS, \Cref{prop:net-qfs} would give a forward-Cauchy net \((h_i)\) in \(E_V\) with Yoneda limit \(\id_V\) and \(\fs_V(h_i)\to1\). The formula forces \(h_i=c_0\) eventually, so its Yoneda limit is \(c_0\ne\id_V\), a contradiction. Thus \(V\) is continuous but not quantitative FS. Together with \Cref{ex:product-dual-fs}, this also shows that quantitative FS is not preserved by taking opposites.
\end{example}

\subsection{Finite separation constructs approximation}\label{sec:approximation}
\begin{theorem}\label{thm:finite-separation-estimate}
Let \(X\) be Yoneda complete, let \(f\colon X\to X\) be continuous, let \(\varphi\in\I X\), and put \(c=\supop_X\varphi\). Then, for every \(x\in X\),
\begin{equation}\label{eq:finite-separation-estimate}
 \fs_X(f)\tn X(x,c)\leq\varphi(f(x)).
\end{equation}
Consequently,
\begin{equation}\label{eq:finite-separation-way-below}
 \fs_X(f)\leq\bigwedge_{x\in X}w_X(f(x),x).
\end{equation}
\end{theorem}

\begin{proof}
Fix \(M\in\Finp(X)\), and put \(u=s_M(f)\) and
\(\psi=f_!\varphi\). Direct images preserve ideals, so
\(\psi\in\I X\); continuity of \(f\) gives \(\supop_X\psi=f(c)\).
For \(z,x\in X\), the inequalities \(u\leq K_M(fx,x)\) and
transitivity give
\[
 u\tn X(z,fx)
 \leq X(z,fx)\tn K_M(fx,x)
 \leq K_M(z,x).
\]
Hence, for every \(z\in X\),
\begin{align*}
 (u\odot\psi)(z)
 &=\bigvee_{x\in X}\varphi(x)\tn u\tn X(z,fx)\\
 &\leq\bigvee_{x\in X}\varphi(x)\tn K_M(z,x)\\
 &=\bigvee_{m\in M}X(z,m)\tn
       \bigvee_{x\in X}\varphi(x)\tn X(m,x)\\
 &=\bigvee_{m\in M}\varphi(m)\tn X(z,m), 
\end{align*}
where the last equality follows from the fuzzy-lower-set law and the diagonal term. The ideal identities and \eqref{eq:fuzzy-supremum} now imply
\begin{align*}
 u
 &=\subdeg_X(\psi,u\odot\psi)\\
 &\leq\subdeg_X(\psi,\bigvee_{m\in M}\varphi(m)\tn X(-,m))\\
 &=\bigvee_{m\in M}\varphi(m)\tn\subdeg_X(\psi,\downarrow_X m) \\
 &=\bigvee_{m\in M}\varphi(m)\tn X(fc,m)
 \leq\varphi(fc).
\end{align*}
For every \(x\in X\), fuzzy order preservation of \(f\) and lowerness of \(\varphi\) give
\[
 u\tn X(x,c)\leq\varphi(fc)\tn X(fx,fc)\leq\varphi(fx).
\]
Joining over \(M\) yields \(\fs_X(f)\tn X(x,c)\leq\varphi(fx)\).
Since \(\varphi\in\I X\) was arbitrary, residuation gives
\[
\fs_X(f)
 \leq\bigwedge_{x\in X}\bigwedge_{\varphi\in\I X}
       \bigl(X(x,\supop_X\varphi)\Rightarrow\varphi(fx)\bigr)
 =\bigwedge_{x\in X}w_X(fx,x).
\]
\end{proof}

\begin{theorem}\label{thm:qfs-continuous}
Every quantitative FS-domain is a continuous quantitative domain. More precisely, if \(\eta_X\) is a finite-separation approximate identity,
\begin{equation}\label{eq:way-below-from-approximate-identity}
 w_X(z,x)=\bigvee_{f\in E_X}\eta_X(f)\tn X(z,fx).
\end{equation}
For each \(x\), this formula defines an ideal \(w_X(-,x)\) with supremum \(x\).
\end{theorem}

\begin{proof}
For \(x\in X\), put
\[
 \rho_x=(\ev_x)_!\eta_X,
 \qquad
 \rho_x(z)=\bigvee_{f\in E_X}\eta_X(f)\tn X(z,fx).
\]
Direct images preserve ideals, and evaluation at \(x\) is continuous.
Thus
\[
 \rho_x\in\I X,\qquad
 \supop_X\rho_x
 =\ev_x(\supop_{E_X}\eta_X)
 =\ev_x(\id_X)=x.
\]
Moreover, by \Cref{thm:finite-separation-estimate}, for every \(\varphi\in\I X\) and \(x,z\in X\),
\begin{align*}
 \rho_x(z)\tn X(x,\supop_X\varphi) 
     &\le \bigvee_{f\in E_X}\fs_X(f)\tn X(z,fx)\tn X(x,\supop_X\varphi)\\
     &\le \bigvee_{f\in E_X}X(z,fx)\tn\varphi(fx)\\
     &\le \varphi(z).
\end{align*}
Therefore, by \eqref{eq:relative-way-below}, we have
\(
    \rho_x(z)\le w_X(z,x).
\)
% \[
%  X(x,x')\tn\rho_x(z)
%  \leq\bigvee_f\eta_X(f)\tn X(z,fx)\tn X(fx,fx')
%  \leq\rho_{x'}(z).
% \]
% Let \(\varphi\in\I X\) and \(c=\supop_X\varphi\).
% By \Cref{thm:finite-separation-estimate} and \(\eta_X\leq\fs_X\),
% \[
%  \rho_x(z)\tn X(x,c)
%  \leq\bigvee_f X(z,fx)\tn\varphi(fx)
%  \leq\varphi(z).
% \]
% Consequently,
% \[
%  X(x,c)\leq\subdeg_X(\rho_x,\varphi)
%  \leq X(\supop_X\rho_x,\supop_X\varphi)=X(x,c).
% \]
% The same pointwise inequality gives \(\rho_x(z)\leq w_X(z,x)\).
Conversely, taking \(\varphi=\rho_x\) in \eqref{eq:relative-way-below} gives
\[
 w_X(z,x)\leq X(x,\supop_X\rho_x)\Rightarrow\rho_x(z)
 =\rho_x(z).
\]
Hence \(w_X(-,x)=\rho_x\in\I X\) and
\(\supop_Xw_X(-,x)=x\).
\end{proof}

% Thus quantitative FS-domains occupy the intended position
% \[
% \FS_T\subsetneq\CONT_T\subsetneq\YC_T.
% \]
% The first inclusion follows from \Cref{thm:qfs-continuous}.

\subsection{Tensor products and function spaces}\label{sec:function-spaces}

The tensor product and continuous function space are complete by \Cref{thm:ambient-real}. We now show that finite separation is preserved by both constructions. Recall that \(X\) is a retract of \(Y\) if there are continuous maps \(e\colon X\to Y\) and \(r\colon Y\to X\) with \(re=\id_X\). Superscripts on \(K_M\) and \(s_M\) specify the underlying fuzzy ordered set.

\begin{proposition}\label{prop:retracts-tensor}
Quantitative FS-domains are closed under retracts in \(\YC_T\) and finite tensor products.
\end{proposition}

\begin{proof}
Let \(e\colon X\to Y\) and \(r\colon Y\to X\) be continuous, with
\(re=\id_X\), and let \(\eta_Y\) be a finite-separation approximate
identity on \(Y\). For \(R(f)=rfe\) and \(M\in\Finp(Y)\),
\begin{align*}
 K^X_{r[M]}(R(f)x,x)
 &\geq\bigvee_{m\in M}X(rfex,rm)\tn X(rm,x)\\
 &\geq\bigvee_{m\in M}Y(fex,m)\tn Y(m,ex) \\
 &\geq s^Y_M(f).
\end{align*}
Thus \(\fs_X(R(f))\geq\fs_Y(f)\). Precomposition and postcomposition
make \(R\colon E_Y\to E_X\) continuous. Set \(\eta_X=R_!\eta_Y\).
Then
\[
 \supop_{E_X}\eta_X=R(\id_Y)=\id_X,
\]
while, for \(h\in E_X\),
\[
 \eta_X(h)=\bigvee_f\eta_Y(f)\tn E_X(h,R(f))
 \leq\bigvee_f\fs_X(R(f))\tn E_X(h,R(f))\leq\fs_X(h).
\]
This proves retract closure.

Let \(X,Y\) be quantitative FS-domains with approximate identities
\(\eta_X,\eta_Y\), and put \(Z=X\otimes Y\).
For \(f\in E_X\) and \(g\in E_Y\), \Cref{thm:ambient-real} gives
\(f\otimes g\in E_Z\) and completeness of \(Z\).
For \(M\in\Finp(X)\) and \(N\in\Finp(Y)\),
\[
 K^Z_{M\times N}((fx,gy),(x,y))
 =K^X_M(fx,x)\tn K^Y_N(gy,y)
 \geq s^X_M(f)\tn s^Y_N(g).
\]
Taking the meet over \(x,y\) and the joins over \(M,N\) gives
\begin{equation}\label{eq:real-tensor-fs-bound}
 \fs_Z(f\otimes g)\geq\fs_X(f)\tn\fs_Y(g).
\end{equation}
Put \(\alpha=\eta_X\otimes\eta_Y\in\I(E_X\otimes E_Y)\), and define
\(\Gamma\colon E_X\otimes E_Y\to E_Z\) by \(\Gamma(f,g)=f\otimes g\). The pointwise fuzzy order satisfies
\[
 E_X(f,f')\tn E_Y(g,g')
 \leq E_Z(\Gamma(f,g),\Gamma(f',g')),
\]
so \(\eta_Z=\Gamma_!\alpha\in\I E_Z\). By \eqref{eq:iterated-fubini} and pointwise ideal suprema,
\begin{align*}
 (\supop_{E_Z}\eta_Z)(x,y)
 &=\supop_\alpha\bigl((f,g)\mapsto(f(x),g(y))\bigr)\\
 &=\supop_{\eta_X}\left(f\mapsto\supop_{\eta_Y}(g\mapsto(f(x),g(y)))\right)\\
 &=\supop_{\eta_X}(f\mapsto(f(x),y))=(x,y).
\end{align*}
Furthermore, \(\alpha(f,g)\leq\fs_Z(\Gamma(f,g))\), and lowerness gives
\[
 \eta_Z(h)=\bigvee_{f,g}\alpha(f,g)\tn E_Z(h,\Gamma(f,g))\leq\fs_Z(h).
\]
Thus \(Z\) is quantitative FS. The singleton has
\(\fs(\id)=1\), so its principal identity ideal supplies the nullary case.
\end{proof}

Jung's proof of closure of FS-domains under function spaces uses finitely many order patterns \cite[Theorem~3]{Jung1990}. In the quantitative setting, finitely many points may give rise to infinitely many patterns of hom-values. The following estimate combines a finite partition of the corresponding finite-dimensional cube with uniform continuity of \(T\) to obtain a uniform finite-separation bound on the function space.

\begin{proposition}\label{prop:real-function-fs}
Let \(X,Y\) be Yoneda-complete fuzzy ordered sets, put \(H=[X,Y]_{\I}\), and let \(f_i\in E_X\), \(g_i\in E_Y\) for \(i=1,2\). Define \(C_i\colon H\to H\) by
\(
 C_i(h)=g_i\circ h\circ f_i.
\)
Then
\begin{equation}\label{eq:real-function-fs}
 \fs_H(C_1C_2)\geq
 \fs_X(f_1)\tn\fs_X(f_2)\tn\fs_Y(g_1)\tn\fs_Y(g_2).
\end{equation}
\end{proposition}
\begin{proof}
For \(h,k\in H\), fuzzy order preservation of \(g_i\) gives
\[
 H(C_i h,C_i k)
 =\bigwedge_xY(g_i h f_i x,g_i k f_i x)\geq H(h,k).
\]
Thus each \(C_i\) is a fuzzy-order-preserving endomap of \(H\). If the right-hand side of \eqref{eq:real-function-fs} is zero, the inequality is immediate. Otherwise, fix
\(0\leq u_i<\fs_X(f_i)\) and \(0\leq v_i<\fs_Y(g_i)\).
Enlarging finite sets by union, choose common \(M\in\Finp(X)\) and
\(N\in\Finp(Y)\) such that, for \(i=1,2\),
\begin{align}
 u_i<s_M(f_i)&\leq\bigvee_{a\in M}X(f_i x,a)\tn X(a,x),
 \label{eq:lower_u}\\
 v_i<s_N(g_i)&\leq\bigvee_{b\in N}Y(g_i y,b)\tn Y(b,y).
 \label{eq:lower_v}
\end{align}

We claim that there exists a finite nonempty set \(L\subseteq H\) such that, for every \(h\in H\), some \(k\in L\) satisfies
\[
H(C_i h,k)\ge u_i\tn v_i,\qquad
H(C_i k,h)\ge u_i\tn v_i
\quad(i=1,2).
\]
Put \(\epsilon=\tfrac12\min_i(s_N(g_i)-v_i)>0\).
Uniform continuity of \(T\) supplies \(\delta>0\) such that
\[
 |q-q'|<\delta\ \Longrightarrow\
 |p\tn q-p\tn q'|<\epsilon,
 \qquad \forall p,q,q'\in[0,1].
\]
Partition \([0,1]^{M\times N}\) into finitely many cells of \(\ell^\infty\)-diameter less than \(\delta\), and choose one map from each cell met by \(\sigma(h)=(Y(b,h(a)))_{a\in M,b\in N}\). The resulting finite set \(L\subseteq H\) satisfies
\[
 \forall h\in H\ \exists k\in L\ \forall a\in M\ \forall b\in N:
 |Y(b,h(a))-Y(b,k(a))|<\delta.
\]
Fix such \(h,k\). For \((h',k')\in\{(h,k),(k,h)\}\), \(a\in M\), and \(i=1,2\),
\begin{align*}
 Y(g_i k'(a),h'(a))
 &\geq\max_{b\in N}Y(g_i k'(a),b)\tn Y(b,h'(a))\\
 &\geq\max_{b\in N}Y(g_i k'(a),b)\tn Y(b,k'(a))-\epsilon\\
 &\geq s_N(g_i)-\epsilon>v_i.
\end{align*}
For each \(i\in\{1,2\}\) and \(x\in X\), choose \(a\in M\) attaining the maximum in
\eqref{eq:lower_u}. Then
\begin{align*}
 u_i\tn v_i
 &\leq X(f_i x,a)\tn Y(g_i k'(a),h'(a))\tn X(a,x)\\
 &\leq Y(g_i k'f_i x,g_i k'(a))\tn Y(g_i k'(a),h'(a))\tn Y(h'(a),h'(x))\\
 &\leq Y(C_i k'(x),h'(x)).
\end{align*}
Hence \(H(C_i k,h)\geq u_i\tn v_i\) and \(H(C_i h,k)\geq u_i\tn v_i\). With \(t=u_1\tn u_2\tn v_1\tn v_2\),
\[
 t\leq H(C_2h,k)\tn H(C_1k,h)
 \leq H(C_1C_2h,C_1k)\tn H(C_1k,h).
\]
The finite subset \(C_1[L]\subseteq H\) therefore gives \(s_{C_1[L]}(C_1C_2)\geq t\). Taking the join over the four thresholds, and using preservation of joins by \(T\), proves
\eqref{eq:real-function-fs}.
\end{proof}

\begin{theorem}\label{thm:pointwise-hom-qfs}
If \(X\) and \(Y\) are quantitative FS-domains, then the function space
\(
  H=[X,Y]_{\I}
\)
is a quantitative FS-domain.
\end{theorem}

\begin{proof}
Let $H=[X,Y]_{\mathcal I}$, and let $\eta_X$ and $\eta_Y$ be finite-separation approximate identities. Put
\[
D=E_X\otimes E_X\otimes E_Y\otimes E_Y,
\qquad
\alpha=\eta_X\otimes\eta_X\otimes\eta_Y\otimes\eta_Y
\in\mathcal I D.
\]
By \eqref{eq:sup-coordinatewise}, 
\[
\sup_D\alpha
=(\mathrm{id}_X,\mathrm{id}_X,\mathrm{id}_Y,\mathrm{id}_Y).
\]
Define $\Gamma:D\to E_H$ by
\[
\Gamma(f_1,f_2,g_1,g_2)(h)
=g_1\circ g_2\circ h\circ f_2\circ f_1.
\]
For Yoneda-complete fuzzy ordered sets \(A,B,C\), \Cref{thm:ambient-real} gives the internal composition
\[
 [B,C]_{\I}\otimes[A,B]_{\I}\longrightarrow[A,C]_{\I},
 \qquad (g,f)\longmapsto g\circ f,
\]
in \(\YC_T\), obtained by currying the composite of evaluation maps. Iterating this construction gives a morphism
\[
 D\otimes H\longrightarrow H,\qquad
 ((f_1,f_2,g_1,g_2),h)
 \longmapsto g_1\circ g_2\circ h\circ f_2\circ f_1.
\]
Its transpose is \(\Gamma\), so \(\Gamma\) is Yoneda continuous.

Put $\theta=\Gamma_!\alpha\in\mathcal I E_H$. By continuity of $\Gamma$,
\[
    \sup_{E_H}\theta
    =\Gamma(\sup_D\alpha)
    =\Gamma(\mathrm{id}_X,\mathrm{id}_X,
             \mathrm{id}_Y,\mathrm{id}_Y)
    =\mathrm{id}_H.
\]
For \(t=(f_1,f_2,g_1,g_2)\), the inequalities \(\eta_X\leq\fs_X\), \(\eta_Y\leq\fs_Y\), and \Cref{prop:real-function-fs} give
\[
\begin{aligned}
 \alpha(t)
 &=\eta_X(f_1)\tn\eta_X(f_2)\tn
   \eta_Y(g_1)\tn\eta_Y(g_2)\\
 &\leq\fs_X(f_1)\tn\fs_X(f_2)\tn
       \fs_Y(g_1)\tn\fs_Y(g_2)\\
 &\leq\fs_H(\Gamma(t)).
\end{aligned}
\]
Equivalently,
\(
\alpha\le\Gamma^*\mathrm{fs}_H.
\)
Since $\mathrm{fs}_H$ is a fuzzy lower set by \Cref{lem:fs-basic}, the direct-image adjunction \eqref{eq:direct-image-adjunction} yields
\(
\theta=\Gamma_!\alpha\le\mathrm{fs}_H.
\)
Therefore $\theta$ is a finite-separation approximate identity on $H$.
\end{proof}

\subsection{Symmetric monoidal closedness}

The closure results now have a categorical consequence. A continuous map of two variables on \(Z\otimes X\) can be regarded as a continuous family, indexed by \(Z\), of continuous maps \(X\to Y\). The two descriptions are related by \(b(z,x)\leftrightarrow\widehat b(z)(x)\) and are inverse to each other.

\begin{theorem}\label{thm:smc-tnorm}
For every continuous \(t\)-norm \(T\), \(\FS_T\) is a symmetric monoidal closed full subcategory of \(\YC_T\) consisting of continuous objects. Its tensor is \eqref{eq:tensor-hom}, its unit is the singleton, and its internal hom is \([X,Y]_{\I}\). The same holds for \(\FS_{T,\bot}\) with all continuous maps.
\end{theorem}
\begin{proof}
By \Cref{thm:qfs-continuous,prop:retracts-tensor,thm:pointwise-hom-qfs},
\[
 \FS_T\subseteq\CONT_T,
 \qquad X\otimes Y\in\FS_T,
 \qquad [X,Y]_{\I}\in\FS_T,
 \quad \forall X,Y\in\FS_T,
\]
and the singleton belongs to \(\FS_T\). Fullness and
\Cref{thm:ambient-real} therefore give the natural bijection
\[
 \FS_T(Z\otimes X,Y)\cong\FS_T(Z,[X,Y]_{\I}),
 \qquad b\longmapsto\bigl(z\mapsto(x\mapsto b(z,x))\bigr).
\]
Evaluation, associativity, symmetry and the unit isomorphisms restrict
from \(\YC_T\).
For pointed \(X,Y\),
\[
 (X\otimes Y)((\bot_X,\bot_Y),(x,y))=1,
 \qquad [X,Y]_{\I}(c_{\bot_Y},h)=1.
\]
The constant map \(c_{\bot_Y}\) is continuous, so both objects are pointed.
The same natural bijection restricts to \(\FS_{T,\bot}\).
\end{proof}

\section{Cartesian closedness}\label{sec:cartesian}
The preceding symmetric monoidal closed structure is available for every continuous \(t\)-norm. Since its tensor need not coincide with the Cartesian product, Cartesian closedness is a further question. This section proves that it holds exactly for the G\"odel \(t\)-norm.

\subsection{Finite Cartesian products}

Since \(\FS_T\) is a full subcategory of \(\YC_T\), closure under Cartesian products implies that these products are also categorical products in \(\FS_T\).

\begin{proposition}\label{prop:cartesian-products}
Quantitative FS-domains are closed under finite Cartesian products.
\end{proposition}

\begin{proof}
Let \(X,Y\in\FS_T\) have approximate identities \(\eta_X,\eta_Y\), and put \(Z=X\times Y\). By \Cref{prop:pointwise-completeness}, \(Z\) and \(E_Z\) are Yoneda complete, and \(f\times g\in E_Z\) for \(f\in E_X\), \(g\in E_Y\). Since \(a\tn b\leq a\wedge b\),
\[
 a\tn b\tn c\tn d=(a\tn c)\tn(b\tn d)\leq(a\wedge c)\tn(b\wedge d).
\]
For \(M\in\Finp(X)\) and \(N\in\Finp(Y)\),
\begin{align*}
 K^Z_{M\times N}((f(x),g(y)),(x,y))
 &\geq\bigvee_{m\in M,n\in N}
 X(f(x),m)\tn X(m,x)\tn Y(g(y),n)\tn Y(n,y)\\
 &=K^X_M(f(x),x)\tn K^Y_N(g(y),y).
\end{align*}
Taking meets over \(x,y\) and joins over \(M,N\) yields
\begin{equation}\label{eq:cartesian-product-fs-bound}
 \fs_Z(f\times g)\geq\fs_X(f)\tn\fs_Y(g).
\end{equation}
Let \(\alpha=\eta_X\otimes\eta_Y\in\I(E_X\otimes E_Y)\), define \(\Gamma\colon E_X\otimes E_Y\to E_Z\) by \(\Gamma(f,g)=f\times g\), and put \(\eta_Z=\Gamma_!\alpha\). The pointwise fuzzy order gives
\[
 E_X(f,f')\tn E_Y(g,g')
 \leq E_Z(f\times g,f'\times g'),
\]
so \(\eta_Z\in\I E_Z\). By \Cref{lem:external-ideal-suprema},
\begin{align*}
 (\supop_{E_Z}\eta_Z)(x,y)
 &=\supop_\alpha((f,g)\mapsto(f(x),g(y)))\\
 &=\supop_{\eta_X}\left(f\mapsto\supop_{\eta_Y}(g\mapsto(f(x),g(y)))\right)\\
 &=\supop_{\eta_X}(f\mapsto(f(x),y))=(x,y).
\end{align*}
Furthermore, \eqref{eq:cartesian-product-fs-bound} and the lower-set property of \(\fs_Z\) give
\[
 \eta_Z(h)
 =\bigvee_{f,g}\alpha(f,g)\tn E_Z(h,f\times g)
 \leq\fs_Z(h).
\]
Thus \(Z\) is quantitative FS. The singleton is the terminal quantitative FS-domain.
\end{proof}

\subsection{A finite obstruction to Cartesian closedness}
The following finite test uses only the defining universal property of a Cartesian exponential. We adapt the finite-test argument of \cite[Section~3]{LaiLuo2026} to the quantitative FS categories, using pointed finite objects and the diagonal scalar obstruction.

\begin{theorem}\label{thm:cartesian-necessity}
Let \(T\) be a continuous \(t\)-norm. If either \(\FS_T\) or \(\FS_{T,\bot}\) is Cartesian closed, then
\[
 p^2=p^3,
 \qquad \forall p\in[0,1].
\]
\end{theorem}

\begin{proof}
Fix \(p\in[0,1]\). Since \(p^3\leq p^2\), it suffices to prove the reverse inequality.

Let \(A=\{\bot,a,b\}\) have fuzzy order defined by
\[
 A(\bot,x)=A(x,x)=1 \quad (x\in A),
 \qquad A(a,b)=p^2,
\]
with all remaining hom-values equal to \(0\). Let \(P=\{\bot,0,1\}\) be defined analogously, with \(P(0,1)=p\). Both are finite separated fuzzy ordered sets. Let
\[
 B=\{0,1,p,p^2,p^3\}\subseteq[0,1]
\]
with repeated elements removed and fuzzy order \(B(r,s)=r\Rightarrow s\). This is a finite separated fuzzy ordered set with least element \(0\).

Define \(f_i\colon A\to B\), for \(i=0,1,2\), by
\[
 f_i(\bot)=0,
 \qquad
 f_i(a)=p^i,
 \qquad
 f_i(b)=p^{i+1}.
\]
The only case requiring verification is the pair \((a,b)\), for which
\[
 p^2\tn f_i(a)=p^{i+2}\leq p^{i+1}=f_i(b),
\]
so each \(f_i\) is fuzzy-order-preserving.

Define \(F,G\colon P\times A\to B\) by the slices
\[
 F(\bot,-)=0,
 \quad F(0,-)=f_0,
 \quad F(1,-)=f_1,
\]
\[
 G(\bot,-)=0,
 \quad G(0,-)=f_1,
 \quad G(1,-)=f_2.
\]
It remains to verify fuzzy order preservation for pairs of the form \(((0,x),(1,x'))\), since the other cases follow from the slice maps and the definitions at the bottom elements. For \(x=x'\),
\[
 p\tn f_i(x)=f_{i+1}(x),
 \qquad \forall i\in\{0,1\},
\]
and for \((x,x')=(a,b)\),
\[
 (p\wedge p^2)\tn f_i(a)=p^2\tn p^i=p^{i+2}=f_{i+1}(b).
\]
Hence \(F\) and \(G\) are fuzzy-order-preserving. All objects involved are finite, so they are quantitative FS and all displayed maps are continuous by \Cref{prop:finite-objects}.

Assume that the relevant quantitative FS category is Cartesian closed, and let \(H=B^A\). Maps from the singleton identify the points of \(H\) with the continuous maps \(A\to B\). Transposing \(F\) and \(G\) gives continuous maps \(\widehat F,\widehat G\colon P\to H\) with
\[
 \widehat F(0)=f_0,
 \quad \widehat F(1)=f_1,
 \qquad
 \widehat G(0)=f_1,
 \quad \widehat G(1)=f_2.
\]
Therefore
\[
 p=P(0,1)\leq H(f_0,f_1),
 \qquad
 p=P(0,1)\leq H(f_1,f_2),
\]
and transitivity in \(H\) yields
\[
 p^2\leq H(f_0,f_2).
\]
Since evaluation \(H\times A\to B\) is fuzzy-order-preserving,
\[
 H(f_0,f_2)\wedge A(a,b)
 \leq B(f_0(a),f_2(b))
 =B(1,p^3)=p^3.
\]
Both factors on the left are at least \(p^2\), so \(p^2\leq p^3\). Every test object has a least element, so the same argument applies to \(\FS_{T,\bot}\).
\end{proof}

\subsection{Classification for continuous \texorpdfstring{\(t\)}{t}-norms}

\begin{theorem}\label{thm:goedel-classification}
Let \(T\) be a continuous \(t\)-norm. The following are equivalent:
\begin{enumerate}[label=\textup{(\roman*)}]
 \item \(\FS_T\) is Cartesian closed;
 \item \(\FS_{T,\bot}\) is Cartesian closed;
 \item \(p\tn q=\min\{p,q\}\) for all \(p,q\in[0,1]\).
\end{enumerate}
\end{theorem}

\begin{proof}
Either \textup{(i)} or \textup{(ii)} implies \(a^2=a^3\) for every \(a\) by \Cref{thm:cartesian-necessity}. Hence
\[
 a^4=a\tn a^3=a\tn a^2=a^3=a^2,
\]
so the continuous map \(s(a)=a^2\) is idempotent. Since \(s(0)=0\) and \(s(1)=1\), the intermediate value theorem gives \(s([0,1])=[0,1]\). For \(r=s(a)\),
\[
 r^2=s(s(a))=s(a)=r.
\]
Thus every scalar is idempotent. If \(p\leq q\), then
\[
 p=p\tn p\leq p\tn q\leq p\tn1=p,
\]
so \(p\tn q=p=\min\{p,q\}\). This proves \textup{(iii)}.

Conversely, under \textup{(iii)} the tensor product is the Cartesian product. The symmetric monoidal closed structure of \Cref{thm:smc-tnorm} is therefore Cartesian closed, proving \textup{(i)} and \textup{(ii)}.
\end{proof}

Thus the product and \L{}ukasiewicz \(t\)-norms yield symmetric monoidal closed quantitative FS categories, but not Cartesian closed ones. For the G\"odel \(t\)-norm, \(\FS_T\) is a Cartesian closed full subcategory of \(\CONT_T\).

\section{Symmetric quantitative FS-domains and compactness}\label{sec:metric}
Finite separation and compactness are linked through classical uniform approximation. Jung and S\"underhauf identify uniformly approximated spaces with FS-domains via specialization order and Scott topology \cite[Theorem~21]{JungSunderhauf1998}. Kummetz characterizes pointed FS-domains by compact approximating \(F\)-posets and identifies the induced topology with the Lawson topology \cite[Theorem~4.3 and Corollary~4.6]{Kummetz1999}.

Here finite separation is defined directly on the enriched structure. We ask what it means when that structure is symmetric. The approximate identity then becomes principal, and its finite-separation condition becomes total boundedness. Combined with Yoneda completeness, this gives compactness of the natural ball topology.

\subsection{Finite separation and total boundedness}
Throughout this section, \(X\) is a nonempty separated symmetric fuzzy ordered set over a continuous \(t\)-norm \(T\). For \(0\leq r<1\), put
\begin{equation}\label{eq:symmetric-open-balls}
 B_r(x)=\{y\in X:X(x,y)>r\},
\end{equation}
and let \(\tau_X\) be the topology generated by these balls.

\begin{lemma}\label{lem:symmetric-fs-cover}
The fuzzy ordered set \(X\) is quantitative FS if and only if it is Yoneda complete and
\begin{equation}\label{eq:symmetric-finite-cover}
 \forall r\in[0,1)\ \exists M\in\Finp(X):
 \qquad X=\bigcup_{m\in M}B_r(m).
\end{equation}
\end{lemma}

\begin{proof}
First let \(A\) be symmetric and let \(\varphi\in\I A\) have supremum \(c\). If \((a_i)\) generates \(\varphi\), then symmetry and the supremum--Yoneda-limit correspondence give
\[
 \varphi(a)
 =\bigvee_i\bigwedge_{j\geq i}A(a,a_j)
 =\bigvee_i\bigwedge_{j\geq i}A(a_j,a)
 =A(c,a)=A(a,c).
\]
Thus every ideal with a supremum is principal.

If \(X\) is Yoneda complete, then \(E_X=[X,X]_{\I}\) is Yoneda complete by \Cref{prop:pointwise-completeness} and symmetric by \eqref{eq:pointwise-hom}. Hence every ideal on \(E_X\) with supremum \(\id_X\) is \(E_X(-,\id_X)\). Consequently, \(X\) is quantitative FS if and only if
\[
 \fs_X(\id_X)=1.
\]
Indeed, necessity follows by evaluating \(\eta_X\leq\fs_X\) at \(\id_X\), while sufficiency follows from \Cref{lem:fs-basic}:
\[
 E_X(f,\id_X)
 =E_X(f,\id_X)\tn\fs_X(\id_X)
 \leq\fs_X(f).
\]

By symmetry,
\[
 \fs_X(\id_X)
 =\bigvee_{M\in\Finp(X)}\bigwedge_{x\in X}
 \bigvee_{m\in M}X(m,x)^2.
\]
If this value is \(1\), then for each \(r\in[0,1)\) there is \(M\in\Finp(X)\) such that every \(x\) admits \(m\in M\) with \(X(m,x)^2>r\), and hence \(x\in B_r(m)\). Conversely, given \(r\in[0,1)\), choose \(s\in[0,1)\) with \(s^2>r\). For a finite cover with centres \(M\), monotonicity gives \(s_M(\id_X)\geq s^2>r\), hence \(\fs_X(\id_X)>r\). Thus \(\fs_X(\id_X)=1\) is equivalent to \eqref{eq:symmetric-finite-cover}.
\end{proof}

\subsection{Compactness and function spaces}
The preceding lemma reduces quantitative FS to Yoneda completeness and finite ball covers. We now identify these with completeness and total boundedness for a separated uniformity.

For \(0\leq r<1\), set
\[
 U_r=\{(x,y)\in X^2:X(x,y)>r\}.
\]
The family \((U_r)_{0\leq r<1}\) is a separated uniformity base. Indeed, it consists of symmetric relations containing the diagonal, is downward directed, and for each \(r\in[0,1)\) one can choose \(s\in[0,1)\) with \(s^2>r\); transitivity gives \(U_s\circ U_s\subseteq U_r\). Separation of \(X\) implies \(\bigcap_{0\leq r<1}U_r\) is the diagonal. The induced topology is \(\tau_X\): if \(y\in B_r(x)\), choose \(s\in[0,1)\) with \(s\tn X(x,y)>r\); then transitivity gives \(B_s(y)\subseteq B_r(x)\).

\begin{theorem}\label{thm:compact-symmetric}
Let \(T\) be a continuous \(t\)-norm and let \(X\) be a nonempty separated symmetric fuzzy ordered set. Then
\[
 X\text{ is quantitative FS}
 \quad\Longleftrightarrow\quad
 (X,\tau_X)\text{ is compact}.
\]
\end{theorem}

\begin{proof}
A net \((x_i)\) is Cauchy for the uniformity \((U_r)_{0\leq r<1}\) if and only if
\[
 \forall r\in[0,1)\ \exists i_0\ \forall j,k\geq i_0:
 X(x_j,x_k)>r.
\]
Every uniformly Cauchy net is forward Cauchy. Conversely, let \((x_i)\) be forward Cauchy and fix \(r\in[0,1)\). Choose \(s\in[0,1)\) with \(s^2>r\), and choose \(i_0\) such that
\[
 X(x_j,x_k)>s,
 \qquad \forall j\geq i_0,\ \forall k\geq j.
\]
For arbitrary \(j,k\geq i_0\), choose \(\ell\geq j,k\). Symmetry gives \(X(x_\ell,x_k)=X(x_k,x_\ell)>s\), and therefore
\[
 X(x_j,x_k)
 \geq X(x_j,x_\ell)\tn X(x_\ell,x_k)
 \geq s^2>r.
\]
Thus the two notions of Cauchy net coincide.

By \cite[Lemma~8.11 and Proposition~9.2]{Zhang2024}, for a Cauchy net in symmetric \(X\), a point \(c\) is a Yoneda limit if and only if \(x_i\to c\) in \(\tau_X\). Hence Yoneda completeness is equivalent to completeness of this separated uniform space.

Condition \eqref{eq:symmetric-finite-cover} is exactly total boundedness of the same uniformity. By the standard compactness theorem for separated uniform spaces \cite[Chapter~6]{Kelley2017}, completeness and total boundedness are equivalent to compactness. The conclusion follows from \Cref{lem:symmetric-fs-cover}.
\end{proof}

% \Needspace{13\baselineskip}
% Here compactness refers to the ball topology \(\tau_X\), not to the Scott topology of the underlying order. The product \(t\)-norm gives the following metric specialization.

% \begin{corollary}\label{cor:compact-metric}
% Let \((X,d)\) be a nonempty metric space, equipped with the fuzzy order
% \[
%  X(x,y)=e^{-d(x,y)}
% \]
% based on the product \(t\)-norm. Then \(X\) is a quantitative FS-domain if and only if \((X,d)\) is compact.
% \end{corollary}

% \begin{proof}
% The triangle inequality makes \(X\) a separated symmetric fuzzy ordered set. For \(0<r<1\),
% \[
%  B_r(x)=\{y\in X:d(x,y)<-\log r\},
% \]
% so \(\tau_X\) is the metric topology. The conclusion follows from \Cref{thm:compact-symmetric}.
% \end{proof}

% The finite-separation condition thus recovers classical FS approximation on Boolean-valued orders and compactness in the symmetric setting. Its preservation by continuous function spaces yields monoidal closedness, while the finite scalar test isolates the additional requirement for Cartesian closedness.

\begin{corollary}\label{cor:complete-metric}
Let \((X,d)\) be a nonempty metric space. Regard \(X\) as a fuzzy ordered set over the product \(t\)-norm by defining
\[
X(x,y)=e^{-d(x,y)}
,\qquad \forall x,y\in X.
\]
Then \(X\) is a quantitative FS-domain if and only if
\((X,d)\) is compact.
\end{corollary}

\begin{proof}
The triangle inequality gives
\[
X(x,y)\tn X(y,z)
=e^{-d(x,y)-d(y,z)}
\leq e^{-d(x,z)}
=X(x,z).
\]
Moreover, \(X(x,x)=1\), and \(X(x,y)=X(y,x)=1\) implies \(x=y\). Thus \(X\) is a separated symmetric fuzzy ordered set.

For \(0<r<1\), the basic open balls defined by \(X\) satisfy
\[
B_r(x)
=\{y\in X:X(x,y)>r\}
=\{y\in X:d(x,y)<-\log r\}.
\]
As \(-\log r\) ranges over all positive real numbers,
the topology \(\tau_X\) coincides with the metric topology
induced by \(d\). The equivalence now follows from
\Cref{thm:compact-symmetric}.
\end{proof}

\begin{corollary}\label{cor:compact-function-space}
If \(X,Y\) are separated symmetric quantitative FS-domains, then \([X,Y]_{\I}\) is compact in its ball topology. In particular, for nonempty compact metric spaces \((X,d_X)\), \((Y,d_Y)\), the space \(\operatorname{Lip}_1(X,Y)\) of nonexpansive maps is compact in the uniform metric
\[
 d_\infty(f,g)=\sup_{x\in X}d_Y(f(x),g(x)).
\]
\end{corollary}

\begin{proof}
The pointwise fuzzy order on \([X,Y]_{\I}\) is separated and symmetric. It is quantitative FS by \Cref{thm:pointwise-hom-qfs}, so \Cref{thm:compact-symmetric} applies.
For the product \(t\)-norm and the exponential fuzzy orders, fuzzy-order-preserving maps are precisely nonexpansive maps. As in the proof of \Cref{thm:compact-symmetric}, forward-Cauchy nets and Yoneda limits are ordinary Cauchy nets and their limits. Nonexpansive maps preserve these limits, so \([X,Y]_{\I}=\operatorname{Lip}_1(X,Y)\). Since compact \(Y\) is bounded, \(d_\infty\) is finite, and
\[
 [X,Y]_{\I}(f,g)=\bigwedge_{x\in X}e^{-d_Y(f(x),g(x))}
               =e^{-d_\infty(f,g)}.
\]
Thus the ball topology is the uniform metric topology.
\end{proof}

% The metric conclusion is a special case of the classical Arzel\`a--Ascoli theorem.

\begin{remark}
If \(X\) is symmetric and Boolean-valued, then separatedness forces its underlying order to be the discrete order. Indeed,
\[
X(x,y)=1
\quad\Longleftrightarrow\quad
X(y,x)=1
\quad\Longrightarrow\quad x=y.
\]
Consequently, \(B_r(x)=\{x\}\) for every \(0\leq r<1\), so \(\tau_X\) is the discrete topology. In this case, \Cref{thm:compact-symmetric} reduces to
\(
X\text{ is quantitative FS}
\text{ if and only if }
X\text{ is finite}.
\)
\end{remark}

\section*{Declaration of generative AI and AI-assisted technologies in the manuscript preparation process}

During the preparation of this manuscript, the authors used ChatGPT to assist with the review of mathematical
arguments, language refinement, consistency checks, and final submission review. The authors subsequently reviewed and revised the AI-assisted content and accept full responsibility for the published article.

\begingroup
\small
\raggedright
\setlength{\bibsep}{4pt plus 1pt minus 1pt}

\endgroup

\end{document}